%
%

\magnification=1200

\font\titfont=cmr10 at 12 pt
\font\headfont=cmr10 at 12 pt



\def\BBB{\bf}

\overfullrule=0in

\def\boxit#1{\hbox{\vrule
 \vtop{%
  \vbox{\hrule\kern 2pt %
     \hbox{\kern 2pt #1\kern 2pt}}%
   \kern 2pt \hrule }%
  \vrule}}

  \def\harr#1#2{\ \smash{\mathop{\hbox to .3in{\rightarrowfill}}\limits^{\scriptstyle#1}_{\scriptstyle#2}}\ }

 \def\GGG{{{\bf G} \!\!\!\! {\rm l}}\ }

\def\bra#1#2{\langle#1, #2\rangle}
\def\ss{\subset}
\def\half{\hbox{${1\over 2}$}}
\def\smfrac#1#2{\hbox{${#1\over #2}$}}
\def\oa#1{\overrightarrow #1}

\def\log{{\rm log}}
\def\Hess{{\rm Hess}}

\def\tr{{\rm tr}}
\def\max{{\rm max}}

\def\span{{\rm span\,}}

\def\Sym{{\rm Sym}^2}
\def\diag{{\rm diag}}

\def\supp{{\rm supp}}

\def\Core{{\rm Core}}

\def\rn{\bbr^n}

\def\Int{{\rm Int}}

\def\Symn{{\Sym(\rn)}}

\def\Theorem#1{\medskip\noindent {\bf THEOREM \bf #1.}}
\def\Prop#1{\medskip\noindent {\bf Proposition #1.}}
\def\Cor#1{\medskip\noindent {\bf Corollary #1.}}
\def\Lemma#1{\medskip\noindent {\bf Lemma #1.}}
\def\Remark#1{\medskip\noindent {\bf Remark #1.}}

\def\Def#1{\medskip\noindent {\bf Definition #1.}}

\def\Ex#1{\medskip\noindent {\bf Example \bf    #1.}}
\def\Qu#1{\medskip\noindent {\bf Question \bf    #1.}}

\def\pf{\medskip\noindent {\bf Proof.}\ }
\def\qed{\hfill  $\vrule width5pt height5pt depth0pt$}
\def\equdef{\buildrel {\rm def} \over  =}

\def\n{\nabla}
\def\w{\wedge}
\def\wee{\wedge\cdots\wedge}

   \def\cp{{\cal P}}   \def\ca{{\cal A}}

\def\ch{{\cal H}}

\def\cp{{\cal P}}

\def\vf{\varphi}

\def\wh{\widehat}

\def\and{\qquad {\rm and} \qquad}

\def\bbr{{\bf R}}

\def\bbp{{\bf P}}
\def\bbd{{\bf D}}

\def\bbf{{\bf F}}

\def\a{\alpha}
\def\b{\beta}
\def\d{\delta}
\def\e{\epsilon}

\def\l{\lambda}
\def\o{\omega}

\def\D{\Delta}
\def\L{\Lambda}
\def\G{\Gamma}
\def\O{\Omega}

\def\psh{plurisubharmonic }

\def\PH#1{\widehat {#1}}

\def\bo{\partial \Omega}

\def\PSH{{ \rm PSH}}

\def\pconv{$p$-convex}

\def\Symn{\Sym(\rn)}

\def\fa{{\rm\ \  for\ all\ }}

\def\AA{1}
\def\BB{2}
\def\BBB{3}
\def\DD{3}

\def\EE{4}
\def\FF{5}

\def\rx{{\cal R}_p(X)}

\ 
\vskip .3in

\centerline{\titfont P-CONVEXITY, P-PURISUBHARMONICTY   }
\smallskip

\centerline{\titfont AND THE LEVI PROBLEM }
\bigskip

\centerline{\titfont F. Reese Harvey and H. Blaine Lawson, Jr.$^*$}
\vglue .9cm
\smallbreak\footnote{}{ $ {} \sp{ *}{\rm Partially}$  supported by
the N.S.F. }

\vskip .4in
\centerline{\bf ABSTRACT} \medskip
  \font\abstractfont=cmr10 at 10 pt

{{\parindent= .8in\narrower\abstractfont \noindent

Three results in $p$-convex geometry are established.
First is the analogue of the Levi problem in several complex variables,
namely: local $p$-convexity implies global $p$-convexity.
The second asserts that the support of a minimal $p$-dimensional current
is contained in the $p$-hull of the boundary union with the ``core'' of the space.
Lastly, the exteme rays in the convex cone of $p$-positive matrices
are characterized.  This is a basic result with many applications.

}}

\vskip .5in

\centerline{\bf TABLE OF CONTENTS} \bigskip

{{\parindent= .1in\narrower\abstractfont \noindent

\qquad \AA. Introduction.\smallskip

\qquad \BB. Plurisubharmonicity.  \smallskip

\qquad \DD.  Convexity, Boundary Convexity, and Local Convexity.  \smallskip

\qquad \EE. Minimal Varieties and  Hulls.\smallskip

\qquad \FF.   Extreme rays in the Convex Cone $\cp_p(V)$.  \smallskip

}}

\vfill\eject

\centerline{\headfont \AA. Introduction.}
\medskip

On any riemannian $n$-manifold there are intrinsic notions of $p$-plurisubharmonicity
and $p$-convexity for   integers $p$ between 1 and $n$. They interpolate between
convexity ($p=1$)  and subharmonicity   ($p=n$) with  
 $p=n-1$ an important case. They arise naturally
in many situations, and their study goes back to H. Wu [Wu].
The object  of this paper is to prove three new results in $p$-geometry.

The central algebraic idea is that of {\sl $p$-positivity} for a quadratic form $Q$
on a finite-dimensional inner product space $V$. By definition $Q$ is 
$p$-positive if the trace of its restriction to every $p$-dimensional subspace $W\ss V$
satisfies $\tr\left\{Q\bigr|_W\right\} \geq 0$. 
This is equivalent to the condition that $\l_1+\cdots+\l_p\geq0$ where
$\l_1\leq\cdots\leq \l_n$ are the ordered eigenvalues of $Q$.
The set of such $Q$ will be denoted $\cp_p(V)$.
On any riemannian manifold $X$, a function
$u\in C^2(X)$ is {\sl $p$-plurisubharmonic} if its riemannian  hessian is $p$-positive.  
An oriented hypersurface
in $X$ is {\sl $p$-convex} if its second fundamental form is $p$-positive. The Riemann curvature 
$R$ of  $X$ is $p$-{\sl positive} if for each tangent vector $v$, the quadratic form 
 $\bra {R_{v,\cdot} \cdot}{v}$ is $p$-positive.
 
 The smooth $p$-plurisubhamonic functions are ``pluri''-subharmonic in the following sense. 
 
  \medskip
 \noindent
 {\bf Theorem \BB.12.}
{\sl A function $u\in C^2(X)$ is $p$-\psh if and only if 
its restriction to every $p$-dimensional minimal submanifold is 
subharmonic in the induced  metric.}

 \medskip
 
The notion of $p$-plurisubhamonicity can be generalized to arbitrary upper
semi-contin-uous $[-\infty,\infty)$-valued functions using standard viscosity test functions (cf. [CIL], [C]).
For $p=1,n$ this recaptures the classical notions of general convex
and subharmonic functions on a riemannian manifold $X$.  This family of upper semi-continuous
$p$-\psh functions,
denoted $\PSH_p(X)$, has many of the useful properties of subharmonic functions
(see Theorem 2.6 in [HL$_{5}$]). Moreover, the Restriction Theorem 2.12 has a non-trivial extension
to general, upper semi-continuous $p$-\psh functions (see [HL$_{6}$]).

The smooth $p$-plurisubharmonic  functions can be used to
 introduce a notion of $p$-convexity as follows.
Given a compact subset $K\ss X$, define the {\sl $p$-convex hull}  of $K$ to be the
set  $\wh K$ of points $x\in X$ such that 
$$
u(x)\ \leq \ \sup_K u
$$
for all smooth $p$-plurisubharmonic functions $u$ on $X$.  Then $X$ is said to be {\sl $p$-convex}
if 
$$
K\ \ss\ss \ X \qquad\Rightarrow \qquad \wh K\ \ss\ss\ X.
$$
The following result was proven in [HL$_{7}$].

\medskip

\centerline
 {\sl A riemannian  manifold $X$ is $p$-convex  if and only if}
 
 \centerline{\sl $X$ admits a smooth $p$-plurisubharmonic proper exhaustion function.}

\medskip

A domain $\O\ss X$ is said to be {\sl  locally $p$-convex} if each  point $x\in \bo$ has a neighborhood
$U$ such that $\O\cap U$ is $p$-convex. Note that  $p$-convex domains are locally  $p$-convex
(see (3.1)).
The following converse is an analogue of Levi Problem in complex analysis, and is  one of the three
new results of this paper.

\medskip
 \noindent
 {\bf Theorem \DD.7.}
 {\sl  Let $\O\ss\ss\rn$ be a domain with smooth boundary.  If $\O$ is locally $p$-convex, then
 $\O$ is $p$-convex.}

\medskip

There is also a notion of  $p$-convexity for the boundary. Let $II$ denote the second
fundamental form of the boundary $\bo$ with respect to the interior normal.  Then the 
boundary $\bo$ is {\sl $p$-convex} if $II_x$ is $p$-positive at each point $x\in\bo$.

\medskip
 \noindent
 {\bf Theorem \DD.9.}
 {\sl  Let $\O\ss\ss\rn$ be a domain with smooth boundary.    If  $\O$ is locally $p$-convex, then  
 $\bo$ is $p$-convex.}

\medskip\noindent
From Theorem \DD.10 one then concludes that for such domains $\O$, 
\medskip
\centerline
{
$\O$ is $p$-convex \qquad$\iff$\qquad
$\O$ is locally $p$-convex \qquad$\iff$\qquad
$\bo$ is $p$-convex.
}
\medskip

A quadratic form $A$ on an inner product space $V$ is said to be {\bf strictly} $p$-positive
if $\tr\left\{A\bigr|_W\right\} > 0$ for all $p$-planes $W\ss V$.  This gives notions of 
strict $p$-plurisubharmonicity, strict $p$-convexity, etc.  In Section 4 a number of results
concerning strictly $p$-convex domains and strictly $p$-convex boundaries are discussed.  
A key concept
here is that of the {\sl core} of $X$, a subset which governs the existence of 
strictly $p$-plurisubharmonic functions and  proper exhaustions (see Remark \EE.4.)

The core contains all compact $p$-dimensional  minimal submanifolds without boundary in $X$.
This result is extended to include non-compact   minimal submanifolds and currents.
A $p$-dimensional rectifiable current $T\in \rx$ on $X$ is called {\sl minimal}
if the first variation of the mass of $T$ is zero with respect to   deformations supported
away from its boundary   $\partial T$ (see Definition \EE.7).

\Cor{\EE.10} {\sl
Suppose $T\in \rx$ is a minimal current, and let $u$ be any smooth $p$-plurisubharmonic
function which vanishes on a neighborhood of $\supp(\partial T)$. Then
$$
 \tr_{\oa T}\left(\Hess \, u \right) \ \equiv \ 0 \qquad{\rm on \ \ } \supp(T).
$$
If $T=[M]$ is the current associated to a connected $p$-dimension minimal submanifold, 
and if the $p$-\psh\ function $u$ and its  gradient  both   vanish at points of $\partial M$,
then}
$$
 u\bigr|_M \ \equiv\ 0 \ \ \ {\rm or, \ if\ \ } \partial M\ =\emptyset,  \ \ 
 u\bigr|_M \ \equiv\ {\rm constant}.
$$

Our second result is the following.

\Theorem{\EE.11} {\sl
Let $K\ss X$ be a compact subset and suppose $T\in \rx$ is a minimal current
such that $\supp(\partial T)\ss K$.  Then}
$$
\supp(T) \ \ss\ \hat K \cup \Core(X).
$$

This leads to the notion of a {\sl minimal surface hull} of a compact set $K\ss X$,
 namely the union of the supports of all
minimal currents $T\in \rx$ whose boundaries are supported in $K$.  Theorem \EE.11 says that
this hull is contained in $\wh K\cup \Core(X)$.

Much of this discussion carries over to minimal (not necessarily rectifiable) $p$-currents.

Our third new result (see Section \FF) describes the extreme rays in the convex 
cone $\cp_p(V)$, defined for each real number $1\leq p \leq n$ by
$$
\cp_p(\rn ) \ \equdef \ \biggl\{A\in\Symn : \l_1(A) + \cdots + \l_{[p]} (A) + (p-[p])\l_{[p]+1}(A) \ \geq\ 0
\biggr\}
\eqno{(\AA.1)}
$$
where $[p]$ denotes the greatest integer $\leq p$ (cf. Remark \BB.9).  
The endpoint cases can be excluded from the discussion
since $\cp_n(V)$ is a half-space (and hence has no extreme rays) while it is well known that 
the extreme rays in $\cp(V) = \cp_1(V) = \{A\geq0\}$ are
generated by the orthogonal projections onto lines. These rays
remain extreme in $\cp_p(V)$ for $1\leq p < n-1$. 
Theorem \FF.1c states that for $1<p<n$ the only other
extreme rays are generated by the elements of $\Sym(V)$ with one negative eigenvalue
$-(p-1)$  and all other eigenvalues 1.

This technical result is more important than it may seem at first glance.
These generators are exactly, up to a positive scale, the second derivatives
of the {\bf Riesz kernel}  $K_p(X)$, which is defined by:

$$
K_p(X) \ =\ 
\cases
{
\ |x|^{2-p} \qquad  \qquad {\rm if}\ \ 1\leq  p<2  \cr
\ \log|x| \  \qquad  \qquad {\rm if}\ \ p=2, \ {\rm and}  \cr
- {1  \over |x|^{p-2}} \qquad  \qquad {\rm if}\ \ 2 <  p \leq n.  \cr
}
\eqno{(\AA.2)}
$$
\noindent
Consequently, an equivalent formulation of Theorem \FF.1c is that following.

\Theorem{\FF.1a. ($1<p<n$)} {\sl  Suppose  $F \ss \Symn$ is a  convex cone
subequation. The Riesz kernel $K_p$ is $F$-subharmonic
if and only if $\cp_p(\rn)\ss F$.
}
\medskip

This result has many applications. One important reason is that it holds for all
{\bf real} numbers $p$  between $1$ and $n$.
 In addition we note that:
\smallskip
\centerline
{\sl
Many of the results from $p$-convex analysis hold for any real $p$, $1\leq p\leq n$. 
}
\smallskip
\noindent
Specifically, since $\cp_p(\rn)\ss\Symn$ is a convex cone, all the results of 
 [HL$_4$] apply.

Finally we note that the basic notions of p-plurisubharmonicity and p-convexity also make sense 
with the grassmann bundle  G(p,TX)
replaced by a closed subset $\GGG\ss G(p,TX)$. There are surprisingly many results which hold in the general context of a ``$\GGG$-geometry''. They are discussed in a companion paper [HL$_7$].



\vfill\eject

\centerline {\headfont \BB. Plurisubharmonicity.}
\medskip

\centerline {\bf  Euclidean Space.}
\medskip

Suppose $V$ is an $n$-dimensional real inner product space, and fix an integer $p$, with
$1\leq p\leq n$.  Let $\Sym(V)$ denote the space of symmetric endomorphisms
of $V$.  Using the inner product, this space is identified with the space of quadratic forms
on $V$.  The notion of $p$-plurisubharmonicity for a smooth function $u$ on $V$ is defined
by requiring that its hessian (i.e., second derivative $D^2_x u$) belong to a certain subset 
$\cp_p(V) \ss\Sym(V)$. To better understand this subset we offer several (equivalent) definitions. 

\Def{\BB.1}  Suppose $A\in \Sym(V)$.  Then $A\in\cp_p(V)$, or $A$ is {\bf $p$-positive}, if the following
equivalent conditions hold.
$$
\tr_WA\ \geq\ 0\qquad\fa\ \  W\in G(p,V)
\leqno{\qquad (1)}
$$
$$
\l_1(A) + \cdots + \l_p(A)  \ \geq \ 0
\leqno{\qquad (2)}
$$
$$
D_A\ \geq\ 0
\leqno{\qquad (3)}
$$

\noindent
where:

\medskip
 
 (1) $G(p,V)$ denotes the set of $p$-dimensional subspaces of $V$, and for $W\in G(p,V)$, the 
 $W$-{\bf trace} of $A$, denoted $\tr_W A$, is the trace  of the restriction $A\bigr|_W$ of $A$
 to $W$,
 
 \medskip
 
 (2)  $\l_1(A) \leq\cdots\leq \l_n(A)$ are the ordered eigenvalues of $A$, so Condition (2)
 says that the sum of the $p$ smallest eigenvalues is $\geq0$, 

 \medskip
 
 (3) $D_A:\L^pV \to \L^pV$ is the linear action of $A$ as a derivation on the space $\L^pV$ of 
 $p$-vectors, i.e., on simple $p$-vectors one has 
 \smallskip
 
 \centerline{$D_A(v_1\wee v_p) = (Av_1)\wedge v_2\wee v_p +
  v_1\wedge (A v_2) \wee v_p+ v_1\w v_2\wee(Av_p)$.}
 \smallskip
The inner product on $V$ induces an inner product on $\L^pV$, and we have 
$D_A\in   \Sym(\L^pV)$, so the notions of non-negativity, $D_A\geq0$, and 
positive definiteness, $D_A>0$,  make sense for $D_A$.\medskip

The proof that condition (1), (2) and (3) are equivalent will be given below.

\Def{\BB.2.   ($p$-plurisubharmonicity)} A smooth function $u$ defined on an open subset
$X\ss\rn$ is said to be  {\bf $p$-plurisubharmonic}  if $D^2_x u  \in \cp_p(\rn)$
for each point $x\in X$.
\medskip

The next result justifies the terminology.

\Prop{\BB.3} {\sl
A function $u\in C^\infty(X)$ is   $p$-plurisubharmonic if and only if the restriction
$u\bigr|_{W\cap X}$ is subharmonic for all affine $p$-planes $W\ss\rn$.
(Here ``subharmonic'' means that $\D_W\left( u\bigr|_{W\cap X} \right) \geq0$ 
where $\D_W$ is the euclidean Laplacian on the affine subspace $W$).
}

\pf This is obvious from Condition  (2) since with $v=u\bigr|_{W\cap X}$, we have 
$\tr_W D^2 u = \D_W v$ on $W\cap X$.\qed

\Remark{\BB.4} The endpoint cases are classical.
\smallskip
\noindent
{${\bf (p = 1)}$ \bf Convex Functions.}  Note that  $A\in \cp_1 \iff  \l_{\rm min}(A) \geq0 \iff A\geq0$, so that 
$u$ is  $1$-plurisubharmonic $\iff$ $u$ is convex.
\smallskip
\noindent
{${\bf (p = n)}$ \bf Classical Subharmonic Functions.} Note that  $A\in \cp_n \iff \tr A\geq0$, so that 
$u$ is  $n$-plurisubharmonic $\iff$ $\D u\geq0$, i.e., $u$ is classically subharmonic.
\smallskip
Consequently, the simplest new case is when $p=2$ in $\bbr^3$ where $u$ is 2-plurisubharmonic $\iff$ the restriction
of $u$ to each affine plane in $\bbr^3$ is classically subharmonic. One generalization of this case has an interesting
characterization.\smallskip
\noindent
{${\bf (p = n-1)}$}  If $p=n-1$, then $*:\L^1 V \to \L^{n-1}  V$ is an isomorphism.
This induces an isomorphism   $\Sym(\L^{n-1} V) \to \Sym(\L^{1} V)$ sending
$
D_A  \ \mapsto  \ (\tr A)I-A.
$
Therefore $u\in C^\infty(X)$ is 
$n-1$-plurisubharmonic if and only if
$$
(\D u)I-\Hess \,u\ \geq\ 0.
$$

\medskip\noindent
{\bf Note.} (a)  It is obvious from Condition (2) that $\cp_p(V)\ss\cp_{p+1}(V)$, or equivalently, if $u$ is   $p$-plurisubharmonic,
then $u$ is $(p+1)$-plurisubharmonic.  In particular, each $p$-plurisubharmonic  function is classically
subharmonic, and every   convex  function is $p$-\psh for all $p$.

\ \ \ (b) The set  $\cp_p(V)$ is a closed convex cone with vertex at the origin.

\medskip
The proof of the equivalence of Conditions (1), (2) and (3) in Definition \BB.1 requires some
elementary facts.  Note that each $p$-plane $W\ss V$ determines a line $L(W) \ss\L^pV$, namely the line
through $v_1\wee v_p$ where $v_1,..., v_p$ is any basis for $W$.  If $e_1,..., e_n$ is an orthonormal 
basis of $V$, we set $e_I = e_{i_1}\wee e_{i_p}$ for $I=(i_1,...,i_p)$ with $i_1< i_2< \cdots <i_p$.

\Lemma {\BB.5} {\sl
Given $A\in\Sym(V)$, consider $D_A \in\Sym(\L^pV)$. Then we have:
\medskip

(a) For all $W\in G(p,V)$,
$$
\tr_W A\ =\ \tr_{L(W)} D_A.
\eqno{(\BB.1)}
$$

(b) If $A$ has eigenvectors $e_1,...,e_n$ with corresponding  eigenvalues $\l_1,...,\l_n$,
then $D_A$ has eigenvectors $e_I$ with corresponding eigenvalues}
$$
\l_I \ = \l_{i_1}+\cdots\l_{i_p}.
\eqno{(\BB.2)}
$$

\pf
For (a), note that,  if $e_1,...,e_p$ is an orthonormal basis of $W$, then $\tr_{L(W)} D_A = \bra {D_A(e_1\wee e_p)}{ e_1\wee e_p}
= \sum_{j=1}^n \bra{e_1 \wee Ae_j \wee e_p}{e_1\wee e_p}
= \sum_{j=1}^n\bra{Ae_j}{e_j} = \tr_W A $.
For (b), compute $D_A e_I =  \l_I e_I$.\qed

\Cor{\BB.6}  {\sl
Suppose $A\in\Sym(V)$ has ordered eigenvalues $\l_1(A)\leq\cdots\leq \l_n(A)$. Then
$$
\inf_{W\in G(p,W)} \tr_W A\ =\ \l_1(A)+\cdots+\l_p(A) \ =\ \l_{\rm min}(D_A),
\eqno{(\BB.3)}
$$
the smallest eigenvalue of $D_A$.}

\pf
Since $D_A$ has eigenvalues $\l_I$ by part (b), the smallest is $ \l_1(A)+\cdots+\l_p(A)=
\tr_{L(\overline W)} D_A$ where $\overline W = \span\{e_1,...,e_p\}$. 
Now the smallest eigenvalue of $D_A$ equals the infimum of $\tr_L D_A$ over all lines in $\L_pV$,
so in this case it is also the infimum over the restricted set of lines of the form $L(W)$ with
$W\in G(p,V)$.  By part (a) in Lemma \BB.5, this proves (\BB.3).\qed

\medskip
The equivalence of Conditions (1), (2) and (3) in Definition \BB.1 is immediate from Corollary \BB.6.

\Def{\BB.7. ($p$-Harmonic)}  A smooth function $u$ defined on an open subset $X\ss\rn$
is {\bf $p$-harmonic} if $D_x^2u \in\partial \cp_p$ for all $x\in X$, or equivalently if $\l_{\rm min}(D_{D^2_x u}) =
\l_1(D^2_x u)+\cdots+\l_p(D^2_x u)=0$  for all $x\in X$.

\Ex{\BB.8. (Radial Harmonics)}
\smallskip
\noindent
{$\bf (p = 1)$}  The function $|x|$ is 1-harmonic on $\rn-\{0\}$.
\smallskip
\noindent
{$\bf (p = 2)$} The function $\log |x|$ is 2-harmonic on $\rn-\{0\}$.
\smallskip
\noindent
{$\bf (3\leq p \leq n)$}    The function $-{1\over |x|^{p-2}}$ is $p$-harmonic on $\rn-\{0\}$.

\pf  
Given a non-zero vector $x\in\rn$, let $P_x \equiv {1\over |x|^2} x\circ x$ denote 
orthogonal projection onto the line through $x$.  One calculates that:
$$
D^2 |x|  \ =\ \smfrac{1}{|x|} \left( I - P_x \right),
\eqno{(\BB.4)}
$$
$$
D^2 \log |x|  \ =\  \smfrac{1}{|x|^2} \left( I - 2 \, P_x \right),
\eqno{(\BB.5)}
$$
$$
D^2 \left( - \smfrac{1}{|x|^{p-2}}\right)  \ =\   \smfrac{(p-2)}{|x|^p}  \left( I - p \, P_x \right).
\eqno{(\BB.6)}
$$

Note that in all cases the function $u(x)$ defined in Example \BB.8 has second derivative
$D^2 u$,  which is  a positive  scalar multiple of  $H\equiv I- p P_x$, 
and that $H$  has one negative eigenvalue 
 $-(p-1)$ and the other  eigenvalues are 1.   By Lemma \BB.5(b) this
implies that the eigenvalues of $D_H$   are  0 and p,  and
in particular, $\l_{\rm min}(D_H)=0$.  \qed

\Remark{\BB.9. (Non-Integer $p$)}
The subset (subequation) $\cp_p(V)$ can be defined for any real number $p$ between
1 and $n$ in such a way that many of the results in this paper continue to hold 
for non-integer values of $p$.  Let $\bar p =[p]$ denote the greatest integer in $p$.  
Then we define $A\in\Sym(V)$ to be {\bf $p$-positive}, or $A\in \cp_p(V)$, if
$$
\l_1(A) + \cdots + \l_{\bar p}(A) + (p-\bar p) \l_{\bar p +1}\ \geq\ 0,
\eqno{(\BB.7)}
$$
where as before $\l_1(A) \leq \cdots \leq \l_n(A)$ denote the ordered eigenvalues of $A$.
To see that $\cp_p(V)$ is a convex cone, one  shows that it is the polar of the set of
$P_{e_1}+ \cdots +P_{e_{\bar p}}+(p-\bar p) P_{e_{\bar p+1}}$ where $e_1,...,e_n$ are orthonormal.

The motivation for this definition of $\cp_p$ is provided by the next remark and Theorem \FF.1.   
These are the only two other places in this paper where non-integer 
values of $p$ are discussed.  In the other places (such as Definition \BBB.1)
the gaps are left to the reader.

\Remark{\BB.10. (The Riesz Kernel)} The family of functions
defined in Example \BB.8 naturally extends by (\AA.2)  to all real numbers $p$ between $1$ and $n$,
and we have the following.
\Lemma{\BB.11}  {\sl  For each real number $p$ with $1\leq p\leq n$},
$$
K_p(x)\ \ {\sl is\ p\, harmonic\,on\ } \rn-\{0\}\ \ {\sl and }\ \  {\sl  p\, plurisubharmonic\,on\ } \rn.
$$

\pf Up to a positive scalar multiple $D^2_x  K_p$ equals $H= I-p P_x$.  As noted above
$D_H \geq0$ and $\l_{\rm min}(D_H)=0$.\qed

\vskip .3in\centerline{\bf Riemannian Manifolds.}\medskip

Suppose $X$ is an $n$-dimensional riemannian manifold. 
Then the euclidean notions above carry over with $V= T_xX$ and the ordinary
hessian of a smooth  function  replaced by the {\bf riemannian hessian}.  
For $u\in C^2(X)$ this is a well defined section of the bundle $\Sym(TX)$
given on tangent vector fields 
$V,W$ by 
$$
(\Hess \, u)(V,W)\ =\ VWu -(\nabla_V W)u,
\eqno{(\BB.8)}
$$
where $\n$ denotes the Levi-Civita connection.
Acting as a derivation,  it determines a well defined section $D_{\Hess\,u}$ of 
$\Sym(\L^p TX)$ for each $p$, $1\leq p\leq n$.

\Def{\BB.2$'$ ($p$-plurisubharmonicity)}  A smooth function $u$ on $X$ is said to be 
 {\bf $p$-plurisubharmonic} if $\Hess_x u$ is $p$-positive at each point $x\in X$
 (see Definition \BB.1).
\medskip

The appropriate geometric objects for restriction are the $p$-dimensional minimal
(stationary) submanifolds of $X$.   In the euclidean case this enlarges
 the family of affine $p$-planes used in Proposition \BB.3 when $1<p<n$.

\Theorem{\BB.12} {\sl
A function $u\in C^2(X)$ is $p$-plurisubharmonic if and only if the restriction of $u$ to every
$p$-dimensional minimal submanifold is subharmonic.
}

\pf
Suppose $M\subset X$ is any $p$-dimensional submanifold, and let $H_M$ denote its mean curvature
vector field.  Then (see Proposition 2.10 in [HL$_2$])
$$
\D_M\left(u\bigr|_M\right)\ =\ \tr_{{{T}}M} \Hess\, u - H_M u.
\eqno{(\BB.9)}
$$
In particular, if $M$ is minimal, then
$$
\D_M\left(u\bigr|_M\right)\ =\ \tr_{{{T}}M} \Hess \,u.
\eqno{(\BB.10)}
$$

It is an elementary fact  (see Lemma \DD.13)  
that for every point $x\in X$ and every $p$-plane $W\ss T_xX$,
 there exists a minimal  submanifold $M$ with $T_xM = W$. This is enough to conclude 
 Theorem \BB.12 from (\BB.10).\qed

\vfill\eject  
 

\centerline{\headfont   \DD.  Convexity, Boundary Convexity, and Local Convexity }
\smallskip
\centerline{\bf    Riemannian Manifolds }\medskip

 Let $\PSH_p^\infty(X)$  denote the smooth
$p$-plurisubharmonic functions on  a riemannian manifold $X$.

\Def{\DD.1}  Given a compact subset $K\subset X$, the {\bf  \pconv\ hull of $K$} is the set
$$
\PH K \ \equiv\ \{x\in X : u(x)\leq \sup_K u  \ \ {\rm for\ all\ \ } u\in \PSH_p^\infty(X)\}
$$

\Prop{\DD.2}
{\sl
If $M\ss X$ is a compact connected $p$-dimensional minimal submanifold with 
boundary $\partial M\neq \emptyset$, then
$$
M\ \ss\ \widehat{\partial M}.
$$
}
\pf Apply Theorem \BB.12 and the maximum principle for subharmonic functions on $M$.\qed

\noindent
{\bf Definition \DD.3.} We say that $X$ is {\bf \pconv}\  if for all compact sets
$K\subset X$, the hull $\PH K$ is also compact.

\Theorem{\DD.4}  {\sl Suppose $X$ is a  riemannian manifold.
Then:
\medskip

(1) \ \ \  $X$ is $p$-convex \qquad$\iff$

\medskip

(2) \ \ \ $X$ admits a smooth $p$-plurisubharmonic proper exhaustion function.
}

\pf
See Theorem 4.4 in [HL$_7$]  for the proof. 
It is exactly the same proof as the one given for Theorem 4.3 in [HL$_2$].\qed
\medskip

Condition (2) can be weakened to a local condition  at $\infty$
in the one-point compactification $\overline X = X \cup \{\infty\}$. This follows from the next lemma.

\Lemma {\DD.5} {\sl
Suppose that $X-K$ admits a smooth $p$-\psh function $v$ with $\lim_{x\to\infty} v(x) = \infty$ where $K$
is compact. Then $X$ admits a smooth $p$-\psh proper exhaustion function which agrees with $v$ near $\infty$.
}

\pf    
This is a special case of Lemma 4.6 in [HL$_7$].
\qed

\vskip .3in
\centerline{\bf   Euclidean Space. }
\medskip
 
We now show that the $p$-convexity of a compact domain with smooth boundary
  in euclidean space is a local condition on the domain near the boundary.  This result is to some degree
  analogous to the Levi Problem in complex analysis, and is one of the three new results of this paper.
 
  \Def{\DD.6}   A domain $\O\ss \rn$ is {\bf locally $p$-convex} if each point $x\in \bo$
  has a neighborhood $U$ in $\rn$ such that $\O \cap  U$ is $p$-convex.
 
 \medskip
 
 Each ball in $\rn$ is $p$-convex, and the intersection of two $p$-convex domains 
 is again $p$-convex.  Therefore
 $$
 {\rm If}\ \ \O\ \ {\rm is\ } p\!-\!{\rm convex, \  then\ } \O \ {\rm is\ locally\ } p\!-\! {\rm convex}.
 \eqno{(\DD.1)}
 $$ 
 Our main result is the converse.
 
 \Theorem {\DD.7} 
 {\sl
 Suppose that $\O$ is a compact domain with smooth boundary.  If $\O$ is locally 
 $p$-convex, then $\O$ is $p$-convex.
 }
 
 \medskip
 Intermediate between local and global convexity is the notion of boundary convexity.
Suppose now that $\bo$ is smooth.

 Denote by $II=II_{\bo}$ the second fundamental form of the boundary with respect to the 
{\bf  inward pointing } normal $n$.  This is a symmetric bilinear form on each tangent space
 $T_x\bo$ defined by 
 $$
II_{\bo}(v,w) \ =\ - \bra {\n_v n}{w}  \ =\   \bra { n}{\n_vW}  
 $$
where $W$ is any vector field tangent to $\bo$ with $W_x=w$.

 \Def{\DD.8}  {\bf The boundary $\bo$ is  \pconv\ at a point $x$} if $\tr_W \{II_{\bo}\} \geq 0$ for all 
 tangential  $p$-planes $W\ss T_x(\bo)$ at $x$.\medskip

 Theorem \DD.7 is the compilation of the following two results.

 \Theorem {\DD.9} 
 {\sl
If  the domain $\O$ is locally   $p$-convex, then its  boundary $\bo$ is $p$-convex.
 }

 \Theorem {\DD.10} 
 {\sl
If  the boundary $\bo$ is   $p$-convex, then the domain $\O$ is $p$-convex.
 }
 
 \medskip
 
 Before proving these two theorems we make some remarks on boundary convexity.
 
 \Remark{\DD.11. (Local defining functions)} Suppose $\rho$ is a smooth function on a neighborhood
 $B$ of a point $x\in\bo$ with $\bo\cap B =\{\rho=0\}$
 and $\O\cap B=\{\rho<0\}$.  If $d\rho$ is non-zero on $\bo\cap B$, then $\rho $ is called a {\bf local
 defining function for $\bo$}.  It has the property that
 $$
 D^2_x \rho \ =\ |\n \rho(x)| II_{x}
  \eqno{(\DD.2)}
 $$
on $\bo\cap B$.  To see this, suppose that $e$ is a vector field tangent to $\bo$ along $\bo$, and note that
$II(e,e) = \bra n {\n_e e} = - {1\over |\n \rho|} \bra {\n\rho}{\n_e e}$ and
$-\bra {\n\rho}{\n_e e} = - (\n_e e)(\rho) =  e(e\rho) - (\n_e e)(\rho) = (D^2\rho)(e,e)$. As a consequence
we have that $\bo$ is $p$-convex at a point $x$ if and only if 
$$
\tr_W D_x^2\rho \ \geq\ 0 \quad {\rm for\ all\ } p\!-\! {\rm planes\ } W  \ {\rm tangent\ to\ }\bo \ {\rm at\ } x
  \eqno{(\DD.3)}
 $$
 where $\rho$ is a local defining function for $\bo$.  Moreover, (\DD.3) is independent of the choice
 of the  local defining function.

 \Remark{\DD.12. (Principal curvatures)}  Let $\kappa_1\leq \cdots\leq \kappa_{n-1}$ denote the 
 ordered eigenvalues of $II_x$.  Then we have that
 $$
 \bo \ \ {\rm is \ } p\!-\!{\rm convex\  at\ } x \qquad \iff\qquad 
 \kappa_1+\cdots + \kappa_p\ \geq\ 0.
  \eqno{(\DD.4)}
 $$
 
 \pf Apply Corollary \BB.6 to $A\equiv II$ with $V \equiv T_x\bo$.\qed
 
 \medskip

 We now give the proof of Theorem \DD.9,   that local $p$-convexity implies boundary $p$-convexity.
 
 \Lemma {\DD.13}
{ \sl
 If $\bo$ is not $p$-convex at a point $x\in \bo$, then there exists an embedded minimal
  $p$-dimensional submanifold $M$ through the point $x$ with
  }
 $$
M-\{x\} \ \ss\ \O\qquad {\rm in\ a\ neighborhood\ of\ }x.  
\eqno{(\DD.5)}
 $$

 \medskip\noindent
 {\bf Proof of Theorem \DD.9}.
 Assume that $\bo$ is not $p$-convex at a point $x\in\bo$. Let $B$ denote the $\e$-ball about $x$.
 It suffices to show that  $\O\cap B$ is not $p$-convex.  This is done by constructing a ``tin can'' inside
 $B$ using Lemma \DD.13.  We can assume that $M$ is a compact manifold with boundary and $M\ss B$.
 
 Let $M_t \equiv M+t\nu$ denote the translate of $M$ by $t\nu$ where $\nu$ is the outward-pointing
 unit normal to $\bo$ at $x$. Choose $r>0$ sufficiently small that each $M_t\ss \O$ for $-r\leq t<0$.
 Let $K$ denote the ``empty tin can'' consisting of the ``bottom'' $M_{-r}$ and the ``label'' 
 $\bigcup_{-r\leq t\leq0} \partial M_t$.  Then $K$ is a compact subset of $\O\cap B$.  
 Let $\hat K$ be its $p$-convex hull in $\O\cap B$.

Since $\partial M_t \ss K$, Proposition \DD.2 implies  that each $M_t \ss\hat K$ for $-r\leq t <0$.  Since $\hat K$ is closed in $\O\cap B$, this proves that $x$ must 
 be in the $\rn$-closure of $\hat K$, i.e., $\hat K$ is not compact.
 Hence, $\O\cap B$ is not $p$-convex.\qed\medskip
 
 \noindent
 {\bf Proof of Lemma \DD.13.}
 Suppose $\bo$ is not \pconv\ at $x$. Then there is a tangent
$p$-plane $W$ to $\bo$ at $x$ with
$$
\tr_W \{II_{\bo} \}\ <\ 0.
\eqno{(\DD.6)}$$
We may assume that $W$ is the plane spanned by eigenvectors of $II$ with 
the  smallest  eigenvalues.
We can then choose euclidean coordinates $(t_1,...,t_n)$ with respect to an orthonormal basis
$e_1,...,e_n$  so that:
\smallskip

(i)  \ \ \ $x$ corresponds to the origin $0$, 

(ii) \  \  $n= e_{n}$ is the outward pointing normal to $\O$ at $x$.

(iii) \ $e_1,..., e_{n-1}$ are the eigenvectors of $II$ at $x$ with eigenvalues  $\kappa_1\leq\kappa_2\leq \cdots\leq \kappa_{n-1}$

(iv) \ \ $W=\span\{e_1,...,e_p\}$  

\medskip
\noindent
In a neighborhood of 0 our domain can be written as 
$$
\O\ =\ \{t_n < f(t_1,...,t_{n-1}) \}.
$$

In particular, $\rho(t) \equiv t_n-f(t_1,...,t_{n-1})$ is a local defining function for $\bo$ near $0\in\bo$.
By Remark \DD.11, since $(\n\rho)(0)=e_n$ is a unit vector, 
$$
D^2_0\rho \ =\ -D^2_0f \ =\ II_0.
\eqno{(\DD.7)}
$$
Hence $f$ has Taylor expansion
 $$
f(t)\ =\ -\half( \kappa_1 t_1^2 +\cdots +  \kappa_{n-1} t_{n-1}^2) + O(|t|^3).
\eqno{(\DD.8)}
$$
 By setting $c\equiv -{1\over p}(\kappa_1 +\cdots +\kappa_{p})$ we obtain a diagonal matrix
 diag$(\kappa_1+c,...,\kappa_p+c)$ with trace zero.  The hypothesis (\DD.6) is equivalent to $c>0$.

We now restrict attention to the linear subspace
$P \equiv  \span\{e_1,...,e_{p}, e_n\} = W \oplus \bbr e_{n} $, and consider
graphs $\{ t_n = g(t_1,...,t_p)\}$ which are minimal hypersurfaces in $P$ (and therefore in $\bbr^n$).
We  apply the following basic lemma, whose proof is left as an exercise.

\Lemma{\DD.14}  {\sl 
 Given $A\in \Sym(\bbr^p)$ with $\tr A=0$, there exists a real analytic function $g$ 
 defined near the origin with $g(0)=0$, $(\n g)(0)=0$ and $D^2_0g = A$ such that $g$ satisfies the minimal
 surface equation.
 }
 \medskip
 
 We can apply this lemma with $A=  -\diag (\kappa_1+c,...,\kappa_p+c)$ obtaining  a minimal surface
 $M=\{(t,g(t)) \in P= \bbr^{p+1} :  |t|<\eta\} \ss \rn$.  The hypothesis $c>0$ implies that $g(t) < f(t_1,...,t_p,0,...,0)$
 if $0<|t|<\eta$ small. This implies that $M-\{0\} \ss\O$, completing the proof of Lemma \DD.13
 and Theorem \DD.9 as well.\qed
 \medskip
 
 Now we commence with the proof of Theorem \DD.10.  Let $\d(x)$ denote the distance from a point
 $x\in\O$ to the boundary $\bo$. 
By the {\bf $\e$-collar of $\bo$}  we shall mean  the set $\{x\in \O : 0<\d(x)<\e\}$.  
 Theorem \DD.10 is immediate from the next result.
 
 \Prop{\DD.15}\smallskip
 {\sl
 
 (1) If $\bo$ is $p$-convex on a neighborhood of $x_0\in \bo$, then
 $-\log \d(x)$ is $p$-\psh on the intersection of  a neighborhood of $x_0$ in $\rn$
 with an $\e$-collar of $\bo$.
 \smallskip
 
 (2) If $-\log \d(x)$ is $p$-\psh on  an $\e$-collar of $\bo$, then $\O$ is $p$-convex.
}

 \medskip
 \noindent
 {\bf  Summary \DD.16.}   From this proposition and Theorems \DD.9 and \DD.10 we conclude that
 $$
 \eqalign
 {
 \O\ \ {\rm is \ locally \ } p \!-\! {\rm convex\ } \quad   
 &\iff\quad  \bo\ \ {\rm is \ } p \!-\! {\rm convex\ }    \cr
  &\iff\quad -\log \d(x)\ \ {\rm is \ } p \!-\! {\rm plurisubharmonic\ }    \cr
  &\iff\quad \O\ \ {\rm is \ } p \!-\! {\rm convex\ }    \cr
  }
 \eqno{(\DD.9)}
$$

\medskip\noindent
{\bf Proof  of  (1).}  
Let $II$ denote the second fundamental form of the hypersurfaces $\{\delta=\epsilon\}$ for $\epsilon\geq0$, 
and let $n= \nabla \d$ denote the inward-pointing normal.  An arbitrary $p$-plane $V$ at a point can be put in a canonical form with basis 
 $$
( \cos \theta) n+ (\sin \theta )e_1, e_2, ... ,  e_p
 $$
where $n, e_1,...,e_p$ are orthonormal.  Set $W\equiv \span\{ e_1,..., e_p\}$, the {\sl tangential
 part } of $V$.

\Lemma {\DD.17} 
$$
\tr_V \, \Hess (-\log \d)\ =\ {1\over \d} \sin^2 \theta \,\tr_W (II) + {1\over \d^2}\cos^2 \theta
$$
\pf See Remark after Proposition 5.13 in [HL$_2$].  \qed 

\medskip
\noindent
{\bf Note.} This formula holds on any riemannian manifold.

\def\l{\kappa}
If $II$ has eigenvalues $\l_1,...,\l_{n-1}$ at a point $x\in\bo$, then let 
$\l_1(\d),...,\l_{n-1}(\d)$ denote the eigenvalues of $II$ at the point a distance
$\d$ from $x$ along the normal line. A proof of  the following can be found in [GT, \S14.6].

\Lemma{\DD.18}  {\sl For small $\d\geq0$ one has}
$$
\l_j(\d)\ =\ { \l_j  \over  1-\d\l_j },  \qquad j=1,...,n-1.
$$

\Cor{\DD.19}  {\sl  Each $\l_j(\d)$ is strictly increasing if $\l_j\neq 0$ and $ \equiv0$ if $\l_j=0$.
}
\medskip

We now combine Lemma \DD.17 with Corollary \DD.19 to conclude that $-\log \d$ is $p$-plurisubharmonic.

\Remark{\DD.20} Note that each $\bo_\e$, where $\O_\e \equiv \{\d>\e\}$, is strictly \pconv\
and $-\log \d$ is strictly $p$-\psh
if and only if $\bo$ has no $p$-flat points, i.e., points where the nullity of $II_{\bo}$ is $\geq p$.
 
\medskip\noindent
{\bf Proof of  (2).}   By Theorem \DD.4 it  suffices to prove the existence
of a continuous exhaustion function $u:\O\to \bbr^+$ which is smooth 
and $p$-plurisubharmonic outside a compact set in $\O$.
Such a function is given by setting $u(x) = \max\{ -\log \d(x),  -\log(\e/2)\}$.
\qed

 \Remark{\DD.21}  It would be interesting to determine if Theorem \DD.9 remains true
 for all real numbers  $p$ between $1$ and $n$.  Most of the other results of this section
 do extend to all such $p$ by [HL$_4$].

 \def\l{\lambda}

\vskip.3in


\vskip .3in
\centerline{\headfont   \EE.  Minimal Varieties and Hulls. }\medskip

There are several notions  of the {\sl $p$-convex hull} of a set, all of which are
intimately related to minimal currents. We begin by recalling the following.
\bigskip
\centerline{\bf Strict Convexity.}
\medskip

Let $X$ be a riemannian manifold which is connected and non-compact.

\Def{\EE.1}
We say that a function $u\in \PSH_p^\infty(X)$ is {\bf strictly $p$-plurisubharmonic\ at a point $x\in X$}
 if $\Hess_x u \in \Int \cp_p(T_xX)$, i.e.,  if one of the following equivalent conditions holds:
 
 \smallskip
 
 (1)\ \ $\tr_W \Hess_x u >0$ for all $W\in G(p,T_xX)$,
 
 \medskip
 
 (2)\ \  $\l_1(\Hess_x u) +\cdots + \l_p(\Hess_x u) >0$,
 
 \medskip
 
 (3)\ \   $D_{\Hess_x u} >0$,
 
 \medskip
 \noindent
 where $\l_1(A)\leq\l_2(A)\leq\cdots$ denote the ordered eigenvalues of $A$.

 \Def{\EE.2}   The manifold $X$ is called  {\bf strictly \pconv}\ if  it admits a \ proper
exhaustion function $u:X\to \bbr$ which is strictly $p$-plurisubharmonic \ at every point, and it is called 
{\bf strictly \pconv\ at infinity} if  it admits a  proper exhaustion function $u:X\to \bbr$
 which is strictly $p$-plurisubharmonic  outside a compact subset.

\Def{\EE.3}  The {\bf $p$-core of $X$} is defined to be the subset
$$
\Core_p(X) \ \equiv\ \{x\in X:  u\ {\rm is\ not\  strict\  at \ } x \ {\rm for \  all \ } u\in \PSH_p^\infty(X)\}
$$

\vfill\eject
\Remark{\EE.4}  
This concept is useful in conjunction with Definition \EE.2.
\medskip

(1)\  $X$ admits a smooth strictly $p$-plurisubharmonic function \quad$\iff$\quad 
\smallskip
\centerline{$\Core(X)=\emptyset$.}

\medskip

(2)\  $X$ is strictly $p$-convex, i.e., $X$ admits a 
smooth strictly $p$-plurisubharmonic proper 

\quad \ \ exhaustion function \quad$\iff$\quad 
\smallskip
\centerline{$\Core(X)=\emptyset$  and $X$ is $p$-convex.}

\medskip

(3)\  $X$ is strictly $p$-convex at infinity \quad$\iff$\quad 
\smallskip
\centerline{$\Core(X)$
is compact and $X$ is $p$-convex.}

\medskip
Part (1) is a special case of Theorem 4.2 in [HL$_7$];  
Part (2) is a special case of 4.8 in [HL$_7$]; and   
Part (3) is a special case of Theorem 4.11 in [HL$_7$].

We note that when
$X$ admits a strictly $p$-plurisubharmonic proper exhaustion function,
 standard Morse Theory implies that $X$ has the homotopy-type of a complex
of dimension $\leq p-1$ (cf. [S], [Wu]).

\Prop{\EE.5} {\sl Every compact $p$-dimensional minimal submanifold $M$
 without boundary in $X$ is contained in  $\Core_p(X)$.  
If instead the boundary  $\partial M \neq \emptyset$ and $M$ is connected,  then}
$$
M\ \subset\  \PH{\partial M}.
$$
\pf  For the first assertion, apply Theorem \BB.12 and the maximum principle to conclude the restriction of 
any smooth $p$-\psh function to $M$ is constant. The second assertion is Proposition \DD.2.\qed
\medskip

This provides an analogue of the support Lemma 3.2 in [HL$_3$].

\Cor{\EE.6} {\sl Suppose $M\ss X$ is a compact $p$-dimensional minimal submanifold with
possible boundary.  Then}
$$
M\ \ss\ \wh {\partial M} \cup \Core(X).
$$

\medskip
\centerline{\bf Minimal Varieties and their Associated Hulls}\medskip

Now we introduce the {\sl minimal current hull} of a compact
set $K$ in a riemannian manifold $X$, and relate it to the $p$-convex hull $\hat K$.
This second hull will be defined using the group $\rx$ 
of $p$-dimensional rectifiable currents with compact support
in $X$ (cf. [F], [Si], [M], etc.). These creatures enjoy many nice properties.  They can
be usefully considered as  compact oriented $p$-dimensional manifolds with singularities
and integer multiplicities, and readers unfamiliar with the general theory can think of  them
simply as submanifolds.

Of importance here is the following general structure theorem. Associated to each  $T\in\rx$
is a Radon measure $\|T\|$ on $X$ and a $\|T\|$-measurable field of unit $p$-vectors
$\oa T$ such that for any smooth $p$-form $\o$ on $X$,
$$
T(\o)\ =\ \int_X \o(\oa T) d\|T\|.
\eqno{(\EE.1)}
$$
(Recall [deR] that the $p$-currents are the topological dual space to the  space
of smooth $p$-forms.) In particular, every $T\in\rx$ has a finite mass
$$
{\bf M}(T) \ =\ \int_X d\|T\|.
$$
{\bf Example.} When $T$ corresponds to integration over a compact 
oriented submanifold with boundary,  of finite volume
 $M\ss X$, one has $\|T\| = \ch_p\bigr|_M$ ($\ch_p$ = Hausdorff measure), $\oa T_x$ corresponds to the 
oriented tangent plane $T_xM$, and ${\bf M}(T) = \ch_p(M)= $ the riemannian volume of $M$.

\Def{\EE.7} A  current $T\in\rx$ is called {\bf minimal} or {\bf stationary} if for all smooth vector fields
$V$ on $X$ which vanish on a neighborhood of the support of $\partial T$, one has
$$
{d\over dt} {\bf M}\left(  (\vf_t)_* T)\right)\biggr|_{t=0} \ =\ 0,
\eqno{(\EE.2)}
$$
where $\vf_t$ denotes the flow generated by $V$ on a neighborhood of the support of $T$.

\medskip

Each smooth vector field on $X$ defines a smooth bundle map  $\ca^V:TX\to TX$
given on a tangent vector $W$ by
$$
\ca^V(W) \ \equdef \ \n_W V.
\eqno{(\EE.3)}
$$
This determines the derivation $D_{\ca^V} : \L^p TX \to \L^p TX$ as in Section \BB.
Proof of the following can be found in [LS] or [L].

\Theorem {\EE.8. (The First Variational Formula)}
{\sl
Fix $T\in \rx$ and let $V$, $\vf_t$ be as above.  Then}
$$
{d\over dt} {\bf M}\left(  (\vf_t)_* T)\right)\biggr|_{t=0} \ =\ \int_X \bra {D_{\ca^V} \oa T}{\oa T} d\|T\|
 \ =\ \int_X \tr_{\oa T}\left( \ca^V \right) d\|T\|.
\eqno{(\EE.4)}
$$

Suppose now that $V=\n u$ for a smooth function $u$ on $X$. Then
$$
\ca^{\n u} = \Hess\, u,
\eqno{(\EE.5)}
$$
considered as an endomorphism of $TX$. To see this note that $\bra {\ca^{\n u}(W) }{U}
= \bra{\n_W(\n u)}{U} = W\bra {\n u}{U} - \bra {\n u}{\n_W U} = (WU - \n_W U)u
= (\Hess\, u) (W, U)$.  Hence, we have the following.

\Theorem{\EE.9} {\sl
If $V = \n u$, then}
$$
{d\over dt} {\bf M}\left(  (\vf_t)_* T)\right)\biggr|_{t=0}
 \ =\ \int_X \tr_{\oa T}\left(\Hess \, u \right) d\|T\|.
\eqno{(\EE.6)}
$$

\Cor{\EE.10} {\sl
Suppose $T\in \rx$ is a minimal current, and let $u$ be any smooth $p$-plurisubharmonic
function which vanishes on a neighborhood of $\supp(\partial T)$. Then
$$
 \tr_{\oa T}\left(\Hess \, u \right) \ \equiv \ 0 \qquad{\rm on \ \ } \supp(T).
$$
If $T=[M]$ where $M$ is a compact connected minimal submanifold of  dimension $p$, 
and if $u$ is a smooth  $p$-\psh function such that $\nabla u\bigr|_{\partial M}=0$,
then}
$$
u\bigr|_M  =\ {\rm constant}.
$$
 
 \pf 
 The first statement follows directly from (\EE.2), (\EE.6) and the fact that $\tr_W \Hess\, u \geq 0$
 on all tangent $p$ planes $W$.  
 
If  $T=[M]$ for a 
 minimal submanifold $M$,  then
 $\tr_{T_xM}( \Hess_x u) = \D_M(u\bigr|_M)$
where $\D_M$ is the Laplace-Beltrami operator of $M$ in the induced metric
(see Proposition 2.10 in [HL$_2$]).  By the First Variational Formula in the smooth
case (e.g. Theorem 1.1 in [L]) we conclude that 
  $u\bigr|_M$ is harmonic on $M$
with constant boundary values (when $\partial M\neq \emptyset$), and the conclusion follows from the maximum principle. 
\qed

\Theorem{\EE.11} {\sl
Let $K\ss X$ be a compact subset and suppose $T\in \rx$ is a minimal current
such that $\supp(\partial T)\ss K$.  Then}
$$
\supp(T) \ \ss\ \hat K \cup \Core(X).
$$
\pf  Suppose $x\notin \hat K$.  Then by the $p$-plurisubharmonic analogue of Lemma 4.2 in [HL$_2$] there exists a 
smooth non-negative $p$-plurisubharmonic function $u$, which is zero on a neighborhood 
of $K$ and satisfies $u(x)>0$, and furthermore, if $x\notin \Core(X)$, 
then $u$ can be chosen to be {\sl strict}
at $x$.  Therefore, $\tr_{\oa T}(\Hess u) >0$ in some neighborhood $U$ of $x$.
Since $\tr_{\oa T}(\Hess u) \geq0$ everywhere, it follows from  (\EE.6), (\EE.2) and minimality that 
 $\|T\|(U)=0$.   Hence, $x\notin \supp(T)$.\qed

\medskip

This result can be rephrased in terms of a second hull defined as follows.

\Def{\EE.12}  Given a compact subset $K\ss X$,  we define the {\bf minimal $p$-current  hull}
to be the set
$$
{\hat K}_{\rm min} \ =\ \bigcup\, \supp(T)
$$
where the union is taken over all minimal $T\in\rx$ with $\supp(\partial T)\ss K$.

\medskip

Note that ${\hat K}_{\rm min} $ contains all compact minimal oriented $p$-dimensional
submanifolds with boundary in $K$.

\Theorem{\EE.11$'$} 
$$
{\hat K}_{\rm min} \ \ss\ \hat K \cup \Core (X).
$$

\medskip

By Remark \EE.4(1), $X$ supports a global strictly $p$-\psh function
if and only if  $\Core(X)=\emptyset$.  Therefore, ${\hat K}_{\rm min} \ss \hat K$ in this case.
For example,   $|x|^2$ is  such a global function on $\rn$.

\Qu{\EE.13}  Suppose $\G\ss \rn$ is a compact $(p-1)$-dimensional submanifold which
bounds exactly one minimal $p$-current in $\rn$ and that current is an oriented  submanifold $M$.
 How close does $\hat \G$ come to approximating $M$?

\vfill\eject

\centerline{\bf General (Not Necessarily Rectifiable) Minimal Currents}
\medskip

Much of what is said above carries over to general compactly supported currents of finite mass.
These are exactly the currents which can be represented as in (\EE.1) with the provision that
the $\|T\|$-measurable  field $\oa T$ of unit $p$-vectors is no longer required to be simple
$\|T\|$-a.e..  Definition \EE.7 makes sense for such currents, and the first variational formula
$$
{d\over dt} {\bf M}\left(  (\vf_t)_* T)\right)\biggr|_{t=0} \ =\ \int_X \bra {D_{\ca^V} \oa T}{\oa T} d\|T\|
$$
holds.  If $V=\n u$ where  $u\in \PSH_p^\infty(X)$, then by  Definition \BB.1 (3) we know that
$D_{\ca^V} \geq 0$.  Furthermore, at any point where $u$ is strict, we have $D_{\ca^V} > 0$.
The arguments for Corollary \EE.10 and  Theorem \EE.11 give the following.

\Prop{\EE.14} 
{\sl Let $T$ be a minimal $p$-dimensional current of finite mass and compact support on $X$,
and let $u$ be any smooth $p$-plurisubharmonic function which vanishes on a neighborhood of $\supp(\partial T)$.
Then 
$$
\bra {D_{\Hess\, u}\oa T}{\oa T} \ =\ 0 \qquad \|T\|-a.e.
$$
Furthermore,
$$
\supp(T) \ \ss\ \wh{\partial T}  \cup  \Core(X).
$$
Thus the minimal current hull $\wh{K}_{\rm min}$ can be expanded to contain the supports of all minimal 
currents with boundary supported in $K$, and Theorem  \EE.11$'$ remains true.}

\medskip
\noindent
{\bf Examples.}  Minimal non-rectifiable currents abound in geometry.  Any   positive 
$(p,p)$-current on a K\"ahler manifold $X$ is minimal. This observation extends to 
positive $\phi$-currents on any calibrated manifold $(X,\phi)$ (see [HL$_1$]).
Any foliation current whose leaves are minimal $p$-submanifolds is a minimal current.

There are two basic cases of smooth minimal currents which are interesting.
Let $T$ be a smooth $d$-closed current of dimension $n-1$ (degree 1). Then 
$T$ is simply a closed 1-form and can be written locally as $T=df$ for a smooth function $f$.
In a neighborhood of any point where $df\neq 0$, the minimality condition is equivalent to 
the 1-Laplace Equation:
$$
d\left(*{{df} \over {\|df\|}}\right) \ =\ 0.
$$
which says that $*{{df} \over {\|df\|}}$ calibrates the level hypersurfaces of $f$.  In particular, the level
sets of $f$ are minimal varieties.

Let $T$ be a smooth $d$-closed  current of dimension 1.  Then $T$ can be expressed on a compactly supported 
1-form $\a$ as 
$T(\a) = \int_X \a(V) d {\rm vol}_X$ where $V$ is a smooth vector field. Minimality is the condition that
$$
\n_V\left({V \over \|V\|}\right)\ =\ 0,
$$
which means exactly that the (reparameterized) flow lines of $V$ are geodesics in $X$, and the $d$-closed condition is equivalent to
$$
{\rm div} (V) \ =\ 0.
$$


\vfill\eject
\centerline{\headfont  \FF. The Extreme Rays in the Convex Cone $\cp_p(V)$. }
\medskip

Recall the classical fact that the extreme rays in $\cp_1(V) \equiv \{A : A\geq0\}$ are
exactly those generated by the orthogonal projections $P_e$ onto the lines spanned
by unit vectors $e\in\rn$.
The purpose of this section is to describe the extreme rays in $\cp_p(V)$ for
other $p$.  Note that $\cp_n(V)$ can be excluded from the discussion since it is a closed
half-space, and hence has no extreme rays. First we state our result in ways that
are more suitable for the many applications (see [HL$_8$] and  [HL$_{9}$] ). 

\Theorem{\FF.1a. ($1<p<n$)} {\sl  The convex cone $\cp_p(\rn) \ss \Symn$
is the smallest convex cone subequation $F$ with the property that the Riesz
kernel $K_p$ is $F$-subharmonic.
}

The second version requires a definition.

\Def{\FF.2}  The {\bf Riesz characteristic} $p_F$ of a subequation $F\ss\Symn$
is defined to be
$$
p_F\ \equiv\ \{p : I-pP_e \in F \ \ \ {\rm for\ all\ } |e|=1\}.
$$

\Theorem{\FF.1b. ($1<p<n$)} {\sl 
Suppose that $F\ss\Symn$ is a convex cone subequation.  Then
}
$$
\cp_p\ss F \qquad\iff\qquad  p\leq p_F.
$$

Finally we state the result in terms of extreme rays.

\Theorem{\FF.1c. ($1<p<n$)} {\sl   The extreme rays in $\cp_p(V)$  are of two types. They are generated by either 
\medskip

\centerline
{
\rm (1) \ \ $I - p(e\circ e) = P_{e^\perp} -(p-1)P_e$   \qquad or  \qquad (2)\ \ $P_e$
} 
\medskip

\noindent
where $e$ is a  unit vector in $V$.  
If $n-1\leq p<n$,  only case (1)  occurs.}
 \medskip

 \noindent
 {\bf Proof of Theorem \FF.1c.}
Under the action of ${\rm O}_n$ on $\Symn$, the set $\bbd \equiv\rn$ of diagonal 
 matrices form an $n$-dimensional cross-section.  
 For any ${\rm O}_n$-invariant set $F\ss\Symn$, the intersection
 $$
 \bbf\ \equiv \ F\cap \bbd
 \qquad {\rm has\ orbit} \qquad
 O( \bbf) \ =\ F.
 $$
 For a convex cone $F\ss \Symn$, let ${{\cal E}}xt( F)$ denote the union of the extreme rays in $F$
 
 \Lemma{\FF.2} {\sl
 If $F\ss \Symn$ is an  ${\rm O}_n$-invariant convex cone and $ \bbf\equiv F\cap \bbd$, then}
 $$
 {{\cal E}}xt(F)\ \subseteq\ O({{\cal E}}xt(  \bbf)).
 $$
 \pf
 Suppose $A\notin O({{\cal E}}xt(  \bbf))$.  Then $A=gDg^t$ with $g\in {\rm O}_n$ 
implies $D\notin {{\cal E}}xt(  \bbf)$.  Thus $D= \a D_0+\b D_1$ with $\a>0$, $\b>0$,
 $D_0, D_1 \in  \bbf$, but $D_0$ and $D_1$ determine different rays.
  Therefore, $A = \a gD_0g^t + \b  gD_1g^t
 =\a A_0 + \b A_1$, $A_0, A_1 \in  O( \bbf) = F$,  but  $A_0$ and $A_1$ determine different rays, proving that 
 $A\notin  {{\cal E}}xt(F)$.\qed

 \medskip
 
 In particular,  ${{\cal E}}xt(  \cp_p) \ss O({{\cal E}}xt(  \bbp_p))$,  so that
  it remains to compute the extreme rays in $\bbp_p \equiv \cp_p\cap \bbd$.
First note that by definition (see Remark \BB.9) we have
 $$
 \bbp_p\ =\ \cp_p\ =\ \{A = {\rm diag}(\l_1,...,\l_n) : \l_1^{\uparrow} + \cdots +  \l_{\bar p}^{\uparrow}
 +(p-\bar p) \l_{\bar p+1}^{\uparrow} \geq 0 \}
 \eqno{(\FF.1)}
 $$
 where $ \l_1^{\uparrow} \leq  \l_2^{\uparrow}\leq \cdots  \l_n^{\uparrow}$ denotes the 
 rearrangement of the $\l_i$'s into ascending order.
 
 \Lemma{\FF.4}
 {\sl
 If $A\in \bbp_p$ is extreme, then $A$ has at most one strictly negative eigenvalue.
 }
 
 \pf
 Suppose $\l_2^{\uparrow} =\l_2^{\uparrow}(A)<0$. 
 To simplify notation we assume the $\l_i$'s are in ascending order and drop the 
 arrows.  Set $\a = \l_1+\l_2 <0$
 and write $\l \equiv (\l_1 , ... , \l_n )$.  Then 
 $\l = sv+(1-s)w$ where $s= {\l_1\over \a}>0$,  $1-s = {\l_2\over \a}>0$, 
 $v=(\a, 0, \l_3, ... , \l_n)$,
 $w=(0, \a, \l_3, ... , \l_n)$,
 and $v,w \in \bbp_p$. Hence, $A$ is not extreme.\qed
 \medskip
 
 We are now reduced to two cases.
 \medskip\noindent
 {\bf One Negative Eigenvalue:}  By rescaling and permuting we may assume 
 $\l_1 = -1$ and $0\leq \l_2\leq \cdots\leq \l_n$ where $A={\rm diag}(\l_1, \l_2, ... ,\l_n)$.
 Set $B={\rm diag}(0,\l_2 \, ... ,\l_n)$. Then 
 $$
 \l_2 + \cdots + \l_{\bar p} +(p-\bar p)\l_{\bar p+1} \ \geq\ 1.
  \eqno{(\FF.2)}
 $$
A similar argument to the one given in the proof of Lemma \FF.4 applies to show that 
if $B$ is extreme in the set of matrices satisfying (\FF.2), then 
$\l_2 = \cdots = \l_n = \mu$ and equality holds in (\FF.2). Therefore,
$( \bar p - 1)\mu  + (p-\bar p)\mu = 1$, that is, $\mu=1/(p-1)$. 
This proves the following.  If $A\in \bbp_p$ is extreme and has one strictly
negative eigenvalue, then after rescaling $A$ and permuting  coordinates,
$A={\rm diag}(-(p-1), 1,  ... , 1)$.
 
  \medskip\noindent
 {\bf All  Eigenvalues Positive:}  Consider the hyperplane
 $ \l_1 + \cdots + \l_{\bar p} +(p-\bar p)\l_{\bar p+1}=1$ intersected with the positive quadrant
 in $\bbr^{\bar p+1}$ (or $\bbr^{\bar p}$ if $p=\bar p$).  The cone on this set is the positive
 quadrant.  Therefore, the only extreme rays of $\bbp_p$ that could possibly appear from this 
 set are the axis rays.  
 \medskip
 
 This proves that the only possible extreme rays in $\cp_p(V)$ are generated by
 $P_e$ and $I-pP_e$ with $|e|=1$.  By the orthogonal invariance of $\cp_p(V)$
 the ray generated by $I-pP_e$, for one unit vector $e$, is extreme if and only if 
 it is extreme for all unit vectors.  Consequently, If $I-pP_e$ is not extreme for one $e$, 
 then the only possible extreme rays are generated by the rank-one projections $P_e$.
 Now $p<n$ implies $\cp_p(V) \cap \{A : \tr A=1\}$ is compact, so  that the extreme rays
 must generate $\cp_p(V)$.  This forces $\cp_p(V) \ss\cp(V)$ which contradicts $1<p$.
 Summarizing, we have that each $I-pP_e$ generates an extreme ray in $\cp_p(V)$.
 
 It remains to show that the axis rays are extreme in $\cp_p(V)$ if and only if $1<p<n-1$.
 This is left to the reader.\qed
 
 \medskip
 
 To see that version a) of Theorem \FF.1 is equivalent to version c), 
 compute that the second derivative $D^2_x K_p$ is, up to a positive scalar,
 equal to $I-pP_x$.
 The equivalence to version b) is straightforward.

\vskip .3in



\centerline{\bf References}

\vskip .2in

\noindent
\item{[C]}   M. G. Crandall,  {\sl  Viscosity solutions: a primer},  
pp. 1-43 in ``Viscosity Solutions and Applications''  Ed.'s Dolcetta and Lions, 
SLNM {\bf 1660}, Springer Press, New York, 1997.

 \smallskip

\noindent
\item{[CIL]}   M. G. Crandall, H. Ishii and P. L. Lions {\sl
User's guide to viscosity solutions of second order partial differential equations},  
Bull. Amer. Math. Soc. (N. S.) {\bf 27} (1992), 1-67.

 \smallskip

\noindent
\item{[deR]}   G. de Rham,   Vari\'et\'es Diff\'erentiables, 
Hermann, Paris, 1960.

\smallskip

\noindent
\item{[F]}  H. Federer, {Geometric Measure Theory},
Springer, 1969.

\smallskip

\item {[GT]}  D. Gilbarg and N. Trudinger, Elliptic Partial Di?erential Equations of Second Order, 
Springer, 1983.

\smallskip

\item {[HL$_{1}$]} F. R. Harvey and H. B. Lawson, Jr., 
 {\sl Calibrated geometries}, Acta Mathematica 
{\bf 148} (1982), 47-157.

 \smallskip

\item {[HL$_{2}$]} F. R. Harvey and H. B. Lawson, Jr., 
 {\sl  An introduction to potential theory in calibrated geometry}, Amer. J. Math.  {\bf 131} no. 4 (2009), 893-944.  ArXiv:math.0710.3920.

\smallskip

\item {[HL$_{3}$]} F. R. Harvey and H. B. Lawson, Jr., {\sl  Duality of positive currents and plurisubharmonic functions in calibrated geometry},  Amer. J. Math.    {\bf 131} no. 5 (2009), 1211-1240. ArXiv:math.0710.3921.

\smallskip

\item {[HL$_{4}$]} F. R. Harvey and H. B. Lawson, Jr.,
{\sl  Plurisubharmonicity in a general geometric context},  Geometry and Analysis {\bf 1} (2010), 363-401. ArXiv:0804.1316

\smallskip

\item {[HL$_{5}$]} F. R. Harvey and H. B. Lawson, Jr., {\sl  Dirichlet duality and the non-linear Dirichlet problem
on Riemannian Manifolds},  J. Diff. Geom. {\bf 88} no. 3  (2011), 395-482.  ArXiv:0907.1981.

\smallskip


\item {[HL$_{6}$]} F. R. Harvey and H. B. Lawson, Jr.,  {\sl  The restriction theorem for fully nonlinear subequations},  to appear in  {\sl Ann. Inst. Fourier},  ArXiv:0912.5220.
\smallskip

\item {[HL$_{7}$]} F. R. Harvey and H. B. Lawson, Jr.,  {\sl  Geometric plurisubharmonicity and convexity -- an introduction},  Advances in Mathematics {\bf 230} (2012), 2428Ð2456.  ArXiv:1111.3875.
\smallskip

\item {[HL$_{8}$]} F. R. Harvey and H. B. Lawson, Jr.,  {\sl  Removable singularities for nonlinear subequations},  ArXiv:1303.0437.
\smallskip

\item {[HL$_{9}$]} F. R. Harvey and H. B. Lawson, Jr.,  {\sl  Radial subequations, isolated singularities
and tangent functions},  (in preparation).
\smallskip

   \noindent
\item{[L]}    H. B. Lawson, Jr.,   Minimal Varieties in Real and Complex Geometry, Les Presses de 
L'Universite de Montreal,  1974.
\smallskip

   \noindent
\item{[LS]}    H. B. Lawson, Jr. and J. Simons,     {\sl On stable currents and their application to global problems in 
real and complex geometry}, Annals of Mathematics 
{\bf 98} (1973), 427-450.

\smallskip

\item {[M]}  F. Morgan,  Geometric Measure Theory -- a Beginner's Guide, Forth Edition,
Elsevier, Academic Press, Amsterdam, 2009.
\smallskip

\item {[S$_1$]} J.-P. Sha, {\sl  $p$-convex riemannian manifolds},
Invent.  Math.  {\bf 83} (1986), 437-447.

\smallskip

\item {[S$_2$]} J.-P. Sha, {\sl  Handlebodies and $p$-convexity},
J. Diff. Geom.  {\bf 25} (1987), 353-361.

\smallskip

\item {[Si]}  L. Simon,  Lectures on Geometric Measure Theory, 
Proc. of the Center for Math. Analysis {\bf 3},  Australian National University, 1983.

\smallskip

\item {[Wu]}   H. Wu,  {\sl  Manifolds of partially positive curvature},
Indiana Univ. Math. J. {\bf 36} No. 3 (1987), 525-548.

\smallskip

\end